\documentclass[12pt]{article}

\def\bl{\begin{eqnarray}}
\def\el{\end{eqnarray}}
\def\bll{\begin{eqnarray*}}
\def\ell{\end{eqnarray*}}
\usepackage{amssymb}
\usepackage{amsmath,amsthm}
\usepackage{amsfonts}
\usepackage{fleqn}
\usepackage{titlesec}
\usepackage[backref=true,dvipdfm,draft=false,colorlinks=true,
bookmarksopen=true,bookmarks=true,bookmarksnumbered=false,
linkcolor=blue,citecolor=blue,urlcolor=blue,pdfauthor=YunBaoHuang]{hyperref}

\title{Maximal right smooth extension chains}

\small \author{Yun Bao Huang\\
 Department of Mathematics\\
 Hangzhou Normal University\\
 Xiasha Economic Development Area\\
 Hangzhou, Zhejiang 310036, China\\
        huangyunbao@sina.com\\
        huangyunbao@gmail.com}
\date{2011.02.26}
\begin{document}
\makeatletter
\newtheorem{defn}{Definition}
\newtheorem{thm}[defn]{Theorem}
\newtheorem{lem}[defn]{Lemma}
\newtheorem{prop}[defn]{Proposition}
\newtheorem{cor}[defn]{Corollary}
\newtheorem{conj}[defn]{Conjecture}
\newtheorem{exmp}[defn]{Example}
\newtheorem{rem}[defn]{Remark}
\makeatother \maketitle
\titlelabel{\thetitle.\;}
\begin{quote}
{\small {\bf Abstract.} If $w=u\alpha$ for  $\alpha\in
\Sigma=\{1,2\}$ and $u\in \Sigma^*$, then $w$ is said to be a
\textit{ simple right extension }of $u$ and denoted by $u\prec w$.
Let $k$ be a positive integer and $P^k(\varepsilon)$ denote the set
of all $C^\infty$-words of height $k$. Set $u_{1},\,u_{2},\cdots,
u_{m}\in P^{k}(\varepsilon)$, if $u_{1}\prec u_{2}\prec \cdots\prec
u_{m}$ and there is no element $v$ of $P^{k}(\varepsilon)$ such that
$v\prec u_{1}\text{ or } u_{m}\prec v$, then $u_{1}\prec u_{2}\prec
\cdots\prec u_{m}$ is said to be a \textit{ maximal right smooth
extension (MRSE) chains }of height $k$. In this paper, we show that
\textit{MRSE} chains of height $k$ constitutes a partition of smooth
words of height $k$ and give the formula of the number of
\textit{MRSE} chains of height $k$ for each positive integer $k$.
Moreover, since there exist the minimal height $h_1$ and maximal
height $h_2$ of smooth words of length $n$ for each positive integer
$n$, we find that \textit{MRSE} chains of heights $h_1-1$ and
$h_2+1$ are good candidates to be used to establish the lower and
upper bounds of the number of smooth words of length $n$
respectively, the method of which is simpler and more intuitionistic
than the previous ones.

 }
\end{quote}
\newpage

\section{Introduction}
Let $\Sigma=\{1, 2\}$, $\Sigma^{*}$ denotes the free monoid over
$\Sigma$ with $\varepsilon$ as the empty word.  If
$w=w_{1}w_{2}\cdots w_{n}$, $w_{i}\in \Sigma$ for $i=1, 2, \cdots,
n$, then $n$ is called the length of the word $w$ and denoted by
$|w|$. For $i=1,\,2$, let $|w|_{i}$ be the number of $i$ which
occurs in $w$, then $|w|=|w|_{1}+|w|_{2}$.

Given a word $w\in\Sigma^{*}$, a \textit{factor} or \textit{subword}
$u$ of $w$ is a word $u\in \Sigma^{*}$  such that $w=xuy$ for
$x,y\in \Sigma^{*}$, if $x=\varepsilon$, then $u$ is said to be a
\textit{prefix} of $w$. A $\mathit{run}$ or $\mathit{block}$ is a
maximum factor of consecutive identical letters.. The
\textit{complement} of $u=u_{1}u_{2}\cdots$$u_{n}\in\Sigma^{*}$ is
the word $\bar{u}=\bar{u}_{1}\bar{u}_{2}\cdots$$\bar{u}_{n}$, where
$\bar{1}=2, \bar{2}=1$.

The Kolakoski sequence $K$ which Kolakoski introduced in
~\cite{Kolakoski}, is the infinite sequence over the alphabet
$\Sigma$

\[K=
\underbrace{1}_{1}\underbrace{22}_{2}\underbrace{11}_{2}\underbrace{2}_{1}
 \underbrace{1}_{1}\underbrace{22}_{2}\underbrace{1}_{1} \underbrace{22}_{2}\underbrace{11}_{2}
 \underbrace{2}_{1}\underbrace{11}_{2}\underbrace{22}_{2}\underbrace{\cdots}_{\cdots=K}\]
which starts with 1 and equals the sequence defined by its run
lengths.

I would like to thank Prof. Jeffrey O. Shallit for introducing me
the Kolakoski sequence $K$ and raising eight open questions on it in
personal communications (Feb. 15, 1990), the fourth and eighth
problems of them are respectively as follows:

(1) Prove or disprove: $|K_i|_1\thicksim |K_i|_2$, which is almost
equivalent to Keane's question.

(2) Prove or disprove: $|K_i|\thicksim \alpha(3/2)^i$ (This would
imply (1)), where $\alpha$ seems to be about 0.873. Does
$\alpha=(3+\sqrt{5})/6$ ?\\
where $K_0=2$ and define $K_{n+1}$ as the string of $1'$ and $2'$
obtained by using the elements of $K_n$ as replication factors for
the appropriate prefix of the infinite sequence $1212\cdots$.

The intriguing Kolakoski sequence $K$ has received a remarkable
attention ~\cite{Br4,Carpi2,Dekking1,Keane,Sing4,Steacy}. For
exploring two unsolved problems, both wether $K$ is recurrent and
whether $K$ is invariant under complement, raised by Kimberling in
\cite{Kim}, Dekking proposed the notion of
\textit{$C^{\infty}$-word} in \cite{Dekking2}. Chv\'{a}tal in
\cite{Chv} obtained that the letter frequencies of
$C^{\infty}$-words are between 0.499162 and 0.500838.

    We say that a finite word $w\in \Sigma^{*}$ in which neither 111 or
222 occurs is \textit{differentiable}, and its \textit{derivative},
denoted by $D(w)$, is the word whose $j$th symbol equals the length
of the $j$th run of $w$, discarding the first and/or the last run if
it has length one.

If a word $w$ is arbitrarily often differentiable, then $w$ is said
to be a \textit{$C^{\infty}$-word} (or \textit{smooth word}) and the
set of all $C^{\infty}$-word is denoted by $\mathcal{C}^\infty$.

A word $v$ such that $D(v)=w$ is said to be a \textit{primitive} of
\textit{w}. Thus 11, 22, 211, 112, 221, 122, 2112, 1221 are the
primitives of 2. It is easy to see that for any word $w \in
\mathcal{C}^{\infty}$, there are at most 8 primitives and the
difference of lengths of two primitives of $w$ is at most 2.

The \textit{height} of a nonempty smooth word $w$ is the smallest
integer $k$ such that $D^k(w)=\varepsilon$ and the height of the
empty word $\varepsilon$ is zero. We write $ht(w)$ for the height of
$w$. For example, for the smooth word $w=12212212$,
$D^{4}(w)=\varepsilon$, so $ht(w)$=4.
\section{Maximal right smooth extension chains}
Let $\mathcal{N}$ be the set of all positive integers and
$P^k(\varepsilon)$ denote the set of all smooth words of height $k$
for $k\in \mathcal{N}$, then
\begin{eqnarray}
P(\varepsilon) &=& \{1,2,12,21\},\label{e1}\\
P^{2}(\varepsilon)&=&\{121,212,11,22,211,122,112,221,2112,1221,1211,12112,\nonumber\\
{}&&2122,21221,1121,21121,2212,12212\}\label{e2}.
\end{eqnarray}
\begin{defn} Let $w,\,u,\,v\in \Sigma^{*}$ if $w=uv$, then $w$ is said to
be a right extension of $u$. Especially, if $v=\alpha\in \Sigma$,
then $w$ is said to be a simple right extension of $u$, and is
denoted by $u\prec w$.
\end{defn}
\begin{defn} Let $u_{1},\,u_{2},\cdots, u_{m}\in
P^{k}(\varepsilon)$, where $k\in \mathcal{N}$.
\begin{eqnarray}
u_{1}\prec u_{2}\prec \cdots\prec u_{m},\label{e3}
\end{eqnarray}
and there is no element $v$ of $P^{k}(\varepsilon)$ such that
\begin{eqnarray}
v\prec u_{1}\text{ or } u_{m}\prec v,\label{e4}\nonumber
\end{eqnarray}
then (\ref{e3}) is said to be a maximal right smooth extension
(\textit{MRSE}) chain of the height $k$. Moreover, $u_1$ and $u_m$
are respectively called the first and last members of the
\textit{MRSE} chain (\ref{e3}).
\end{defn}

Let $H^{k}$ denote the set of all \textit{MRSE} chains of the height
$k$. For $\xi\in H^{k}$, $\xi=u_{1}\prec u_{2}\prec \cdots\prec
u_{m}$, the complement of $\xi$ is $\bar{u}_{1}\prec
\bar{u}_{2}\prec \cdots\prec \bar{u}_{m}$, and is denoted by
$\bar{\xi}$. It is clear that $\bar{\xi}$ is also a \textit{MRSE}
chain of the height $k$. In addition, for $A\subseteq H^{k}$,
$\bar{A}=\{\bar{\xi}:\;\xi\in A\}$.
\begin{defn} For $\xi\in H^{k+1}$, $\xi=u_{1}\prec u_{2}\prec
\cdots\prec u_{m}$, where $k\in \mathcal{N}$. If there is an element
$\eta=v_{1}\prec v_{2}\prec \cdots\prec v_{n}\in H^{k}$ such that
$u_{1},\,u_{2},\cdots, u_{m}$ are all the primitives of
$v_{1},\,v_{2},\cdots, v_{n}$, then $\xi$ is said to be a primitive
of $\eta$.
\end{defn}
For example, $\xi=121\prec 1211\prec 12112\in H^{2}$ is a primitive
of $\eta=1\prec 12\in H^{1}$, $\bar{\xi}=212\prec 2122\prec 21221\in
H^{2}$.

For a set $A$, let  $|A|$ denote the cardinal number of $A$. Next we
establish the formula of the number of the members of $H^{k}$. For
this reason, let
\begin{eqnarray}
H^{k}_{1}&=&\{\xi\in H^{k}:\;\xi=u_{1}\prec u_{2}\prec \cdots\prec
u_{m}\textit{ and \textit{first}}(u_{1})=1\};\label{e5}\\
H^{k}_{2}&=&\{\xi\in H^{k}:\;\xi=u_{1}\prec u_{2}\prec \cdots\prec
u_{m}\textit{ and \textit{first}}(u_{1})=2\}\label{e6}.
\end{eqnarray}
It immediately follows that
\begin{eqnarray}
H^{k}_{1}&=&\bar{H}^{k}_{2};\label{e7}\\
H^{k}&=&H^{k}_{1}\cup H^{k}_{2};\label{e8}\\
|H^{k}_{1}|&=&|H^{k}_{2}|.\label{e9}
\end{eqnarray}
From (\ref{e1}) and (\ref{e2}) we have
\begin{eqnarray}
H^{1}&=&\{1\prec 12, 2\prec 21\};\label{e12}\\
H^{1}_{1}&=&\{1\prec 12\};\nonumber\\
H^{1}_{2}&=&\{2\prec 21\};\nonumber\\
H^{2}&=&\{121\prec 1211\prec 12112,\; 212\prec 2122\prec 21221,\;
11\prec 112\prec 1121,\;\nonumber\\
{}&&22\prec 221\prec 2212,\; 211\prec 2112\prec 21121,\; 122\prec
1221\prec 12212\}; \label{e13}\\
H^{2}_{1}&=&\{121\prec 1211\prec 12112,\; 11\prec 112\prec 1121,\;
122\prec 1221\prec 12212\};\nonumber\\
H^{2}_{2}&=&\{212\prec 2122\prec 21221,\; 22\prec 221\prec 2212,\;
211\prec 2112\prec 21121\}.\nonumber
\end{eqnarray}
Thus from (\ref{e12}) and (\ref{e13}),  we see that every
\textit{MRSE} chain of height $k$ is uniquely determined  by its
first member $u_{1}$ and each member of $P^{k}(\varepsilon)$ exactly
belongs to one \textit{MRSE} chain of height $k$ for $k=1,2$ and
\begin{equation}
|H^{2}|=3|H^{1}|.\label{e10}
\end{equation}
Actually, the above result holds for every $k\in \mathcal{N}$.\\
\begin{thm}\label{thm3} $H^{k}$ is stated as above. Then each
member of $P^{k}(\varepsilon)$ exactly belongs to one \textit{MRSE}
chain of height $k$, that is, $H^{k}$ gives a partition of the
smooth words of height $k$ and
\begin{equation}
|H^{k}|=2\cdot 3^{k-1}\text{ for all } k\in \mathcal{N}.\label{e11}
\end{equation}
\end{thm}
\noindent{\bf Proof.} We proceed by induction on $k$. From
(\ref{e10}) it follows that (\ref{e11})  holds for $k=1,2$. Assume
that (\ref{e11}) holds for $k=n-1\geq 1$.

Now we consider the case for $k=n$. Since for each $\eta=u_{1}\prec
u_{2}\prec \cdots\prec u_{m}\in H^{n-1}_{1}$, from the definition 2
and (\ref{e5}), we see that
$first(u_{1})=first(u_{2})=\cdots=first(u_{m})=1$, and
$u_{i+1}=u_{i}\alpha$ where $i=1,2,\cdots,m-1$, $\alpha=1,2$. Thus
if $\alpha=1$ then the two primitives $p(u_{i+1})$ of $u_{i+1}$ are
\begin{eqnarray}
p(u_{i+1})&=&\bar{\beta}\Delta^{-1}_{\beta}(u_{i+1})\gamma\nonumber\\
&=&\bar{\beta}\Delta^{-1}_{\beta}(u_{i})\bar{\gamma}\gamma\nonumber\\
&=&p(u_{i})\gamma\nonumber, \textit{ where } \beta,\,\gamma\in
\Sigma,
\end{eqnarray}
so $p(u_{i})\prec p(u_{i+1})$.

If $\alpha=2$ then the four primitives $p_{t}(u_{i+1})$ of $u_{i+1}$
are
\begin{eqnarray}
p_{t}(u_{i+1})&=&\bar{\beta}\Delta^{-1}_{\beta}(u_{i+1})\gamma^{t}\nonumber\\
&=&\bar{\beta}\Delta^{-1}_{\beta}(u_{i})\bar{\gamma}^{2}\gamma^{t}\nonumber\\
&=&p(u_{i})\bar{\gamma}\gamma^{t}\nonumber, \text{ where }
\beta=1,2,\; t=0,1,
\end{eqnarray}
hence $p(u_{i})\prec p_{0}(u_{i+1})\prec  p_{1}(u_{i+1})$.
Therefore, $\eta$ has exactly two primitives and the primitives of
$u_{1},\,u_{2},\cdots, \text{ and }u_{m}$ all occur in the two
primitives of $\eta$.

For example, $\eta=121\prec 1211\prec 12112\in H^2_1$ has exactly
two primitives:

$\mu=121121\prec 1211212\prec 12112122\prec 121121221$ and
$\bar{\mu}$.

Analogously, we can see that each member $\eta$ of $H^{n-1}_{2}$ has
exactly four primitives and the primitives of $u_{1},\,u_{2},\cdots,
\text{ and }u_{m}$ all occur in the four primitives of $\eta$.

For example, $\eta=212\prec 2122\prec 21221\in H^{2}_{2}$ has
exactly four primitives:

$\xi_{1}=22122\prec 221221\prec 2212211\prec 22122112\prec
221221121$;

$\xi_{2}=122122\prec 1221221\prec 12212211\prec 122122112\prec
221221121$ and $\bar{\xi}_{1},\;\bar{\xi}_{2}$.

Thus, by the induction hypothesis, it follows from (\ref{e8}) and
(\ref{e9}) that
\begin{eqnarray}
|H^{n}|&=&|H^{n}_{1}|+|H^{n}_{2}|\nonumber\\
&=&2\cdot |H^{n-1}_{1}| +4\cdot |H^{n-1}_{2}|\nonumber\\
&=&3\cdot (|H^{n-1}_{1}|+|H^{n-1}_{2}|)\nonumber\\
&=&3\cdot |H^{n-1}|\nonumber\\
&=&2\cdot 3^{n-1}.\hspace{.3cm}\Box\nonumber
\end{eqnarray}
\section{The number of smooth words of length $n$}
Let  $\gamma(n)$ denote the number of smooth words of length $n$ and
$p_{K}(n)$  the number of subwords of length $n$ which occur in $K$.

Dekking in~\cite{Dekking2} proved that there is a suitable positive
constant $c$ such that \\$c\cdot n^{2.15}\leq\gamma(n)\leq n^{7.2}$
and brought forward the conjecture that there is a suitable positive
constant $c$ satisfying $p_{K}(n)\sim c\cdot
n^{q}(n\rightarrow\infty)$, where $q=(\log 3)/\log (3/2)$.

Recall from~\cite{Weakley} that a $C^{\infty}$-word $w$ is
\emph{left doubly extendable} (LDE) if both $1w$ and $2w$ are
$C^{\infty}$, and a $C^{\infty}$-word $w$ is \emph{fully extendable}
(FE) if $1w1, 1w2, 2w1$, and $2w2$ all are $C^{\infty}$-words. For
each nonnegative integer $k$, let $A(k)$ be the minimum length and
$B(k)$ the maximum length of an FE word of height $k$.

Weakley in ~\cite{Weakley}  proved that there are positive constants
$c_{1}$ and $c_{2}$ such that for each $n$ satisfying $B(k-1)+1\leq
n \leq A(k)+1$ for some $k$, $c_{1}\cdot n^{q}\leq \gamma(n) \leq
c_{2}\cdot n^{q}$.

It is a pity that we don't know how many positive integers $n$
fulfil the conditions required. Set $\gamma'(n) = \gamma(n+1) -
\gamma(n)$, Weakley in ~\cite{Weakley} gave
\begin{equation}  \label{eeq1}
\gamma(n) = \gamma(0) + \sum_{i=0}^{n-1} \gamma'(i)  \textit{ for }
n \geq 2.
\end{equation}

Let $F(n)$ denote the number  of LDE-words of height $n$, Shen and
Huang in \cite[Proposition 3.2]{Huang4} established
\begin{equation}  \label{eeq2}
F(n)=4\cdot 3^{n-1}\textit{ for each positive integer } n.
\end{equation}

Huang and Weakley in \cite{Huang3} combined (\ref{eeq1}) with
(\ref{eeq2}) to show that
\begin{thm}[\cite{Huang3} Theorem 4]\label{thm1}
Let $\xi$ be a positive real number and $N$ a positive integer such
that for all LDE words $u$ with $|u| > N$ we have $|u|_{2}/|u| >
(1/2) - \xi$. Then there are positive constants $c_{1}, c_{2}$ such
that for all positive integers $n$,
\[
c_{1}\cdot n^{\frac{\log 3}{\log((3/2) + \xi + (2/N))}} < \gamma(n)
< c_{2}\cdot n^{\frac{\log 3}{\log({(3/2) - \xi})}}.
\]
\end{thm}
Let $\gamma_{a,b}(n)$ denote the number of smooth words of length
$n$ over 2-letter alphabet $\{a,\,b\}$ for positive integers $a<b$,
Huang in \cite{Huang5} obtained
\begin{thm}\label{t1}
For any positive real number $\xi$ and positive integer $n_{0}$
satisfying $|u|_b/|u|>\xi$ for every LFE word $u$ with $|u|>n_0$,
there exist two suitable constants $c_{1}\text{ and }\, c_{2}$ such
that
\begin{eqnarray}
c_{1}\cdot  n^{\frac{\log (2b-1)}{\log
(1+(a+b-2)(1-\xi))}}\leq\gamma_{a,b}(n)\leq c_{2}\cdot n^{\frac{\log
(2b-1)}{\log (1+(a+b-2)\xi)}}\nonumber
\end{eqnarray}
for every positive integer $n$.
\end{thm}
Since there are the minimum height $h_1(n)$ and maximum height
$h_2(n)$ of smooth words of length $n$ for each positive integer
$n$, so the lengths of smooth words of the height $h_1(n)-1$ must be
less than $n$ and the length of smooth words of the height
$h_2(n)+1$ must be larger than $n$. Now we are in a position to use
the number of \textit{MRSE} chains of the suitable height $k$ to
bound the number of smooth words of length $n$, which is simpler
than the ones used in \cite{Huang3,Huang5}. The estimates for the
heights of smooth words of length $n$ are borrowed from
\cite{Huang5} in the following proof.
\begin{thm}\label{thm2}  For any positive number $\theta$ and $n_{0}$
satisfying $|u|_2/|u|>\theta$ for $|u|>n_{0}$, there are suitable
positive constant  $c_{1}, c_{2}$ such that
\begin{eqnarray}
c_{1}\cdot n^{\frac{\log 3}{\log (2-\theta)}}\leq\gamma(n)\leq
c_{2}\cdot n^{\frac{\log 3}{\log (1+\theta)}} \textit{ for any
positive integer } n.\nonumber
\end{eqnarray}
\end{thm}

\noindent{\bf Proof. } It is obvious that
\begin{equation}\label{eq4}
|w|\leq |D(w)|+|D(w)|_{2}\textit{ for each smooth word } w.
\end{equation}

First, since $|u|_{2}/|u|>\theta$ for $|u|\geq n_0$, from
(\ref{eq4}) one has
\[|w|\geq (1+\theta)|D(w)| \textit{ for } |D(w)|_{2}/|D(w)|>\theta,\]
which implies
\begin{eqnarray}
|D(w)|<\alpha |w|\textit{ for }|w|\geq N_{0},\label{eq1}
\end{eqnarray}
where $N_{0}$ is a suitable fixed positive integer such that
$|D(w)|\geq n_0$ as soon as $|w|\geq N_{0}$, $\alpha=1/(1+\theta)$.
Since there are finitely many smooth words of length less than
$N_{0}$, from (\ref{eq1}) we see that there exists a suitable
nonnegative integer $l$ such that
\begin{eqnarray}
|D(w)|<\alpha |w|+l\textit{ for each smooth word}\,.\label{eqq1}
\end{eqnarray}

Let $k_0$ be the least integer such that the length of smooth words
of height $k_0$ is larger than $\frac{l}{1-\alpha}$  and $r$ be the
smallest length of smooth words of height $k_0$.  Let $k$ be the
height of the smooth words $w$ such that $ht(w)\geq k_0$, then
$ht(D^{k-k_0}(w))= k_0$. So, from (\ref{eqq1}), we get
\begin{eqnarray}
r&\leq & |D^{k-k_0}(w)|\nonumber\\
&<& \alpha|D^{k-k_0-1}(w)|+l\nonumber\\
&<& \alpha^2|D^{k-k_0-2}(w)|+\alpha l+l\nonumber\\
&\cdots&\nonumber\\
&<&\alpha^{k-k_0}|w|+\alpha^{k-k_0}l+\cdots
+\alpha^2 l +\alpha l+l\nonumber\\
&<&\alpha^{k-k_0}|w|+\frac{l}{1-\alpha}.\nonumber
\end{eqnarray}
Thus
\[
(1/\alpha)^{k-k_0}<\frac{|w|}{\lambda}, \textit{ where
}\lambda=r-\frac{l}{1-\alpha},\] which means
\begin{eqnarray}
ht(w)=k<\frac{\log |w|}{\log (1/\alpha)}+k_0-\frac{\log
\lambda}{\log (1/\alpha)}.\nonumber
\end{eqnarray}
Since there are only finitely many smooth words satisfying $ht(w)<
k_0$, so there is a suitable constant $t_2$ such that
\begin{eqnarray}
ht(w)<\frac{\log |w|}{\log (1/\alpha)}+t_2\textit{ for each smooth
word}\,.\label{eq2}
\end{eqnarray}
Therefore, the maximal height $h_2(n)$ of all smooth words of length
$n$ satisfies
\begin{eqnarray}
h_2(n)&\leq& \frac{\log n}{\log (1+\theta)}+t_{2}.\label{eq3}
\end{eqnarray}

Put $k=h_2(n)+1$, then the length of every smooth word of height $k$
is greater than $n$, so each smooth word of length $n$ can be right
extended to get a \textit{MRSE} chain of height $k$, which suggests
$\gamma(n)\leq |H^k|$. Consequently, from (\ref{e11}) and
(\ref{eq3}) it follows the desired upper bound of $\gamma(n)$.

Second, since the
 complement of any smooth word is a smooth word of the
 same length, the theorem's hypothesis implies that $|D(w)|_{1}/|D(w)|
 \geq \theta$, so $|D(w)|_{2}/|D(w)| \leq 1 - \theta$. From
(\ref{eq4}) it follows that
\begin{equation}\label{e14}
|w|\leq \beta |D(w)|+q \textit{ for each } C^\infty\textit{-word}
\;w,
\end{equation}
where $\beta=2-\theta,\,q$ is a suitable positive constant. Thus
\begin{eqnarray}
|w|\leq\beta^{k-1} |D^{k-1}(w)|+q\frac{\beta^{k-1}-1}{\beta-1}<
2\beta^{k-1}+\frac{q\beta^{k-1}}{\beta-1}
=(2+\frac{q}{\beta-1})\beta^{k-1}=t\beta^{k-1},\nonumber
\end{eqnarray}
where $t=2+q/(\beta-1)$, $k$ is the height of $|w|$. Wherefore, the
length $|w|$ of a smooth word $w$ with height $k$ is less than
$t\beta^{k-1}$ and $k-1>(\log |w|-\log t)/\log \beta$. Hence, the
smallest height $h_1(n)$ of smooth words of length $n$ meets
\begin{eqnarray}
h_1(n)> \frac{\log n}{\log (2-\theta)}+t_1, \textit{ where
}t_1=1-\frac{\log t}{\log \beta} \label{eq31}.
\end{eqnarray}

Then the length of all smooth words with height $m=h_1(n)-1$ is less
than $n$, which means that the length of the last member
$\textit{last }(\xi)$ is less than $n$ for each $\xi\in H^m$. Since
each smooth words of length no more than $n-1$ can be extended right
to  a smooth word of length $n$, we see $\gamma(n)\geq |H^{m}|$.
Herewith, from (\ref{e11}) and (\ref{eq31}) we  get the desired
lower bound of $\gamma(n)$.   $\Box$
\section{Concluding remarks}
Let $a\text{ and }b$ be positive integers of different parities and
$a<b$. Lately, Sing in \cite{Sing4} conjectured:

\textit{There are positive constants $c_1,\,c_2$ such that}
 \begin{eqnarray}
c_{1}\cdot n^\delta\leq\gamma_{a,b}(n)\leq c_{2}\cdot n^\delta,\,
\text{ where }\delta=\frac{\log(a+b)}{\log((a+b)/2)}.\nonumber
\end{eqnarray}
Theorem \ref{t1} means Sing's conjecture should be revised to  be of
the following form
\begin{eqnarray} c_{1}\cdot n^\theta\leq\gamma_{a,b}(n)\leq
c_{2}\cdot n^\theta,\, \text{ where
}\theta=\frac{\log(2b-1)}{\log((a+b)/2)}.\nonumber
\end{eqnarray}

For 2-letter alphabet $\Sigma=\{a,b\}$ with $a<b$, let
$P^j(\varepsilon)$ denote the set of smooth words of height $k$ for
$j\in \mathcal{N}$. For $\alpha\in\Sigma$, set
\[\xi_i=\alpha^i\prec\alpha^i\bar{\alpha}\prec \alpha^i\bar{\alpha}^2\prec \cdots
\prec \alpha^i\bar{\alpha}^{b-1}\text{ for }1\leq i\leq b-1.\] and
\[H^1=\{\eta|\eta=\xi_i\text{ or }\bar{\xi_i},\,i=1,2,\cdots,b-1\}.\]
Let $H^2$ denote the set of primitives of the members in $H^1$, then
it is easy to see $H^2$ constitutes a partition of
$P^2(\varepsilon)$. So continue, we can define the set $H^k$ for
each $k\in \mathcal{N}$ and $H^k$ constitutes a partition of
$P^k(\varepsilon)$. Using the method similar to Theorem \ref{thm2},
we could establish the corresponding result to Theorem \ref{t1}.

\end{document}